\documentclass[11pt]{article}
\usepackage{latexsym}
\usepackage{amsthm}
\usepackage{amssymb}
\usepackage{a4wide}
\def\halmos{\mbox{\,\rule{0.25cm}{0.30cm}}}
\newtheorem{theorem}{Theorem}
\newtheorem{corollary}{Corollary}
\newtheorem{lemma}{Lemma}
\newtheorem{definition}{Definition}
\newtheorem{proposition}{Proposition}
\title{Akivis superalgebras and speciality}
\author{Helena Albuquerque \thanks{Departamento de Matem\'atica,
Faculdade de Ci\^encias e Tecnologia da Universidade de Coimbra,
Apartado 3008, 3001-454 Coimbra, Portugal; E-mail: {\tt
lena@mat.uc.pt}. Financial Support by CMUC-FCT gratefully
acknowledged} \and Ana Paula Santana
\thanks{Departamento de Matem\'atica,
Faculdade de Ci\^encias e Tecnologia da Universidade de Coimbra,
Apartado 3008, 3001-454 Coimbra, Portugal; E-mail: {\tt
aps@mat.uc.pt}. Financial Support by CMUC-FCT gratefully
acknowledged} }
\begin{document}

\maketitle

\hskip2cm{\it In honour of Ivan Shestakov on the occasion of his
60th birthday}

\begin{abstract}
In this paper we define Akivis superalgebra and study enveloping
superalgebras for this class of algebras,  proving an analogous of
the PBW Theorem.

 Lie and Malcev superalgebras are examples of
Akivis superalgebras. For these par\-ti\-cu\-lar superalgebras, we
describe the connection between the classical enveloping
superalgebras and the corresponding generalized concept defined in
this work.
\end{abstract}

{\bf AMS-Classification}: 15A63,  17A70

\section{Introduction}

\begin{definition}
The supervector space $M=M_0\oplus M_1$ is called an Akivis
superalgebra if it is endowed with two operations:\begin{itemize}
\item a bilinear superanticommutative map $\rm{[\- ,  \- ]}$ that
induces on M a structure of superalgebra;\item a trilinear map $A$,
compatible with the gradation (i.e, $A(M_\alpha,  M_\beta,
M_\gamma)\subseteq M_{\alpha+\beta+\gamma}$, all $\alpha,  \beta,
\gamma \in Z_2)$, satisfying the following identity:\newline $[[x,
y], z]+(-1)^{\alpha(\beta+\gamma)} [[y,  z],  x]+
(-1)^{\gamma(\beta+\alpha)} [[z,  x],  y]= A(x,  y,
z)+(-1)^{\alpha(\beta+\gamma)} A(y,  z,  x)+
+(-1)^{\gamma(\beta+\alpha)} A(z,  x,  y)-(-1)^{\alpha \beta} A(y,
x, z)-(-1)^{\alpha(\beta+\gamma)+\beta \gamma} A(z,  y,  x)-
(-1)^{\gamma \beta} A(x,  z,  y),  $ \newline for homogeneous
elements $x\in M_\alpha,   y\in M_\beta,   z\in M_\gamma$,   all
$\alpha, \beta, \gamma \in Z_2$. \end{itemize} This superalgebra
will be denoted in this work by $(M,  [ ,   ],  A)$,  or simply $M$,
if no confusion arises.
\end{definition}

This definition is a generalization of the notion of Akivis algebra
presented by I.    Shestakov in [8].
    In fact,  the even part of an Akivis superalgebra is an Akivis algebra.    Akivis algebras were introduced by M.
        A.    Akivis in [1] as local algebras of local analytic loops.

\vspace{4mm} In this paper, we consider  Akivis superalgebras over a
field $K$ of characteristic different from 2 and 3.    It is our aim
to study the enveloping superalgebra of an Akivis superalgebra and
to prove an analogous of the PBW Theorem. Our approach is similar to
the one used in [8] but,  as it is expected, the superization of the
results imply more elaborated calculations and arguments. This is
particularly evident in the definition of the superalgebra $\tilde
V(M)$ studied in Section 5.

\vspace{5mm}

Given a superalgebra $W$ with multiplication $( ,   )$,   we will
denote by  $W^-$ the superalgebra with underlying supervector space
$W$ and multiplication $[ ,   ]$ given by $[x,  y]=(x,
y)-(-1)^{\alpha\beta}(y,  x)$ for homogeneous elements $x\in
W_\alpha$,   $y\in W_\beta$ (and extended by linearity to every
element of $W$). It is known that if $W$ is an  associative
superalgebra then the superalgebra $W^-$ is a Lie superalgebra,  and
if $W$ is an  alternative superalgebra then $W^-$ is a Malcev
superalgebra. Standard calculations show that if $W$ is any
superalgebra,   the superalgebra $W^-$ is an Akivis superalgebra
for the trilinear map $A(x,  y,  z)= (xy)z-x(yz)$.    This
superalgebra will be denoted in this work by $W^A$.

 We recall that  a superalgebra $S$ is said to be special if it is isomorphic to $U^-$ for
 some  superalgebra $U$.

 It is well known that every  Lie superalgebra is isomorphic to a superalgebra $S^-$, where
 $S$ is an associative superalgebra.    In fact,  for any Lie superalgebra $L$,
 let $T(L)$ denote
 the associative tensor superalgebra  of the vector space $L$, and  consider its bilateral ideal
 $I$ generated by the homogeneous elements $\,x\otimes y-(-1)^{\alpha\beta} y\otimes
 x-[x,y],$ all $
 x\in M_\alpha,   y\in M_\beta,   \alpha,  \beta \in Z_2\,.\,$ Then  the associative superalgebra
 $T(L)/I$ is the universal enveloping superalgebra of $L$ and $L$ is isomorphic to  a
 subsuperalgebra of  $({T(L)/I})^-$.

It is an open problem to know if a  Malcev superalgebra is
isomorphic to  $S^-$ for some alternative superalgebra  $S$.    This
problem was solved only partially in  [3] and [4]. There  it was
shown that, in some cases, a Malcev algebra is isomorphic to a
subalgebra of ${Nat(T)}^-$ for an algebra $T,  $ where ${Nat(T)}$
denotes the generalized alternative nucleous of $T$.

In this work,    we prove that  an  Akivis superalgebra $M$  defined
over a field $K$ of characteristic different from 2 and 3 is
isomorphically  embedded in $\tilde U(M)^A$, where $\tilde U(M)$ is
its  enveloping superalgebra. So $M$ is special. \vspace{5mm}

\vspace{5mm}

\section { Examples of Akivis superalgebras}

Lie superalgebras and more generally Malcev superalgebras are Akivis
superalgebras. For the first class,   we consider the trilinear map
A to be the zero map and, for the second class,  we take $A(x,  y,
z)=1/6\,  SJ(x,  y,  z).    $ Here
 $SJ(x,  y,  z)$ denotes the superjacobian
$$
SJ(x,  y,  z)= [[x,  y],  z]+(-1)^{\alpha(\beta+\gamma)} [[y,  z],
x]+ (-1)^{\gamma(\beta+\alpha)} [[z,  x],  y],
$$
of the homogeneous elements   $x\in M_\alpha,   y\in M_\beta,   z\in
M_\gamma$,  $(\alpha,  \beta,  \gamma \in Z_2).   $

\vspace{5mm} Next,  we give two examples of Akivis superalgebras
which are not included in these classes.

Consider the algebra of octonions $O$ as the algebra obtained by
Cayley-Dickson Process from the quaternions $Q$,  with the $Z_2$
gradation $$O_0=Q= <1,  e_1,  e_2,  e_3> \,\,\,\,\,
\mbox{and}\,\,\,\,\, O_1=e_4Q = < e_4,  e_5,  e_6,  e_7>.   $$  The
multiplication table of the Akivis superalgebra $O^A$ is shown
below.    Note that the even part of this superalgebra is the simple
Lie algebra $sl(2,  K)= <e_1, e_2, e_3>$ together with $1\in Z(O)$.

\begin{center}
\begin{tabular}{|c||c|c|c|c|c|c|c|c|}
  \hline
&1&$e_1$&$e_2$&$e_3$&$e_4$&$e_5$&$e_6$&$e_7$\\  \hline\hline
1    &0&0      &0      &0&0&0&0&0\\
$e_1$&0&0      &$-2e_3$&2$e_2$&$-2e_5$&$2e_4$&$2e_7$&$-2e_6$\\
$e_2$&0&$2e_3$ &0      &$-2e_1$&$-2e_6$&$-2e_7$&$2e_4$&$2e_5$\\
$e_3$&0&$-2e_2$&$2e_1$&$0$&$-2e_7$&$2e_6$&$-2e_5$&$2e_4$\\
$e_4$&0&$2e_5$ &$2e_6$ &$2e_7$&$-2$&0&0&0\\
$e_5$&0&$-2e_4$&$2e_7$ &$-2e_6$&0&-2&0&0\\
$e_6$&0&$-2e_7$&$-2e_4$&$2e_5$&0&0&-2&0\\
$e_7$&0&$2e_6$ &$-2e_5$&$-2e_4$&0&0&0&-2\\
\hline
\end{tabular}\end{center}

\vskip5mm

This superalgebra is an Akivis superalgebra that is neither a Lie
nor a Malcev superalgebra.    Indeed,  we have that $SJ(e_3,  e_7,
e_2)\not=0$ and $((e_4 e_2)e_3)e_5- ((e_2 e_3)e_5)e_4\not= (e_4 e_3)
(e_2e_5)$. \vspace{2mm}

The second class of examples can be obtained using antiassociative
superalgebras (nonassociative $Z_2$-graded quasialgebras). Consider
$D$ a division algebra and $n,  m$ natural numbers.    In [3] the
authors studied the superalgebra  $\widetilde{Mat}_{n,  m}(D)$  of
the $(n+m)\times(n+m)$ matrices over $D$,   with the chess-board
$Z_2$-grading
\begin{eqnarray*}
\widetilde{Mat}_{n,  m}(D)_0 =\left\{
                                \left(\begin{array}{cc}
                                  a & 0 \\
                                  0 & b \\
                                \end{array}\right)
 : a\in Mat_n(D), \,   b\in Mat_m(D)\right\},\\
\widetilde{Mat}_{n,  m}(D)_1 =\left\{ \left(\begin{array}{cc}
                                  0 & v \\
                                  w & 0 \\
                                \end{array}\right)
 : v\in Mat_{n\times m}(D), \,  w\in Mat_{m\times n}(D)\right\}
\end{eqnarray*}
and with multiplication given by
\[
\pmatrix{a_1&v_1\cr w_1&b_1\cr}\cdot \pmatrix{a_2&v_2\cr w_2&b_2\cr}
= \pmatrix{a_1a_2+v_1w_2&a_1v_2+v_1b_2\cr
w_1a_2+b_1w_2&-w_1v_2+b_1b_2\cr} \ .
\]
This superalgebra is antiassociative,  with even part isomorphic to
$Mat_n(D)\times Mat_m(D)$.\newline ${\widetilde{Mat}_{n,  m}(D)}^A$
is an Akivis superalgebra that is  neither a Lie  nor a Malcev
superalgebra.    In fact,  if $E_{i,  j}$ denotes the $(n+m)\times
(n+m)$ matrix with $(ij)$ entry equal to 1 and all the other entries
equal to 0,   then  $SJ(E_{1,   n+1},  E_{n+1,  1},  E_{1,
n+1})\not=0$ and $2((E_{1,   n+1}E_{n+1,  1})E_{1,   n+1})E_{n+1,
1}-$ \newline $-((E_{n+1,  1}E_{1,   n+1})E_{n+1,  1})E_{1,
n+1}\not= E_{1,   n+1}^2 E_{ n+1,  1}^2$.

In the case $m=n=1$,   $\tilde Mat_{1,  1} (D)^A$ has abelian even
part and multiplication given by the following table,   where
$a=E_{11},  b= E_{11}-E_{22},  x=E_{12},  y=E_{21}$: \vskip1cm

\begin{center}\begin{tabular}{|c||c|c|c|c|}
  \hline
&a&$b$&$x$&$y$\\  \hline\hline
$a$&0&0      &x&$-y$\\
b    &0&0      &2x      &-2y\\
$x$&-x&-2x  &0&$b$\\
$y$&y&2y&$b$ &0      \\
\hline
\end{tabular}\end{center}

\vspace{5mm}

\section  {Enveloping superalgebra of an Akivis superalgebra}

In this section we construct and study the universal enveloping
superalgebra of an Akivis superalgebra.

Given  the  Akivis superalgebras $(M,   [ ,   ],  A)$ and $(N,    [
,   ]',  A'),  $ by an Akivis homomorphism we mean  a superalgebra
homomorphism of degree 0, $f:M\rightarrow N$,  such that, for all
$x,  y,  z\in M, \,f(A(x,  y,  z))= A'(f(x),  f(y),  f(z)).$

\begin{definition}
Let $M$ be an Akivis superalgebra.    A pair $(\tilde U,  \iota)$ is
an universal enveloping  superalgebra of $M$ if $\tilde U$ is a
superalgebra and $\iota:M\rightarrow {\tilde U}^A$ is an Akivis
homomorphism satisfying the following condition: given any
superalgebra $W$ and any Akivis homomorphism $\theta:M\rightarrow
W^A$,   there is a unique superalgebra homomorphism of degree 0,
$\tilde \theta:\tilde U\rightarrow W$, such that $
\theta=\tilde\theta \iota$.
\end{definition}

In a similar way to the one used in the classical case for Lie
superalgebras we can prove the following:

\begin{proposition}
1) The universal enveloping  superalgebra of an Akivis superalgebra
is unique up to isomorphism;

2) Let $M$ be an Akivis superalgebra and $(\tilde U,  \iota)$ its
universal enveloping  superalgebra.    Then the  superalgebra
$\tilde U$ is generated by $\iota(M)$ and $K$;

3)Consider two Akivis superalgebras $M_1$ and $M_2$ with universal
enveloping  superalgebras $(\tilde U_1,  \iota_1)$ and $(\tilde U_2,
\iota _2), $ respectively.    If there is an Akivis homomorphism
$\phi: M_1\rightarrow M_2$ then there is a superalgebra homomorphism
$\tilde \phi:\tilde U_1\rightarrow \tilde U_2$ such that $\tilde
\phi\iota_1=\iota_2\phi$.

4)Let $M$ be an Akivis superalgebra with  universal enveloping
superalgebra $(\tilde U,  \iota)$.    Let $I$ be a graded ideal in
$M$ and let $J$ be the graded ideal of $\tilde U$ generated by
$\iota(I)$.    If $m\in M$,   the map $\lambda:m+I \rightarrow
\iota(m)+J$ is an Akivis homomorphism of $M/I$ in $(\tilde U/J)^A$
and $(\tilde U/J,  \lambda)$ is the universal enveloping
superalgebra of $M/I$.
\end{proposition}

\vspace{5mm}

Next we will construct the universal enveloping  superalgebra of the
Akivis superalgebra $(M,   [ ,   ],  A)$.    We start by considering
the nonassociative $Z$-graded tensor algebra of $M$
$$
\tilde T(M)= \oplus_{n\in Z} \tilde T^n(M),
$$
 where $\tilde T^n(M)=0$,  if $n<0$,   $\tilde T^0(M)= K$,
   $\tilde T^1(M)=M$ and  $\tilde T^n(M)= \oplus_{i=1} ^{n-1} \tilde T^i(M)\otimes \tilde T^{n-i}(M),$
   for $n\geq 2, $ with multiplication  defined by
$xy=x\otimes y\,$.   $\,\tilde T(M)$ is a superalgebra with the
$Z_2$-gradation   $$\tilde T(M)=(\oplus_{n\ge 0} \tilde
T^n(M)_0)\oplus (\oplus_{n\ge 1} \tilde T^n(M)_1), $$
 where $$\tilde T^n(M)_\gamma = \oplus _{i=1}^{n-1}\oplus _{\alpha+\beta=\gamma} (\tilde T^i(M)_\alpha \otimes \tilde T^{n-i}(M)_\beta),
   \,\,\gamma\in Z_2.$$
Let $I$ be the $Z_2$-graded ideal of $\tilde T(M)$ generated by the
homogeneous elements $$x\otimes y -(-1)^{\alpha\beta} y\otimes x-[x,
y] \,\,\,\,\mbox{and}\,\,\,\, (x\otimes y)\otimes z -x\otimes
(y\otimes z)-A(x, y, z),$$ for $x\in M_\alpha,   y\in M_\beta, z\in
M_\gamma$. The quocient algebra $\tilde U(M)=\tilde T(M)/I$ is a
superalgebra with the natural $Z_2$-gradation induced by the graded
ideal $I$. Consider the map $\iota:M\rightarrow \tilde
T(M)\rightarrow \tilde U(M)$ obtained by the composition of the
canonical injection with the quocient map.    It is obvious that
$\iota$ is an Akivis homomorphism between $M$ and $\tilde U(M)^A$.

\begin{proposition}
The superalgebra $(\tilde T(M)/I,  \iota)$ is the universal
enveloping  superalgebra of the Akivis superalgebra $M$.
\end{proposition}

\noindent {\bf Proof:} Given any  superalgebra $W$ and   an  Akivis
homomorphism $f:M\rightarrow W^A$,  we need to prove that there is
a unique homomorphism $\tilde f:\tilde U(M)\rightarrow W$ such that
$f=\tilde f \iota$.

Using the universal property of tensor products it is easy to see
that there is a unique superalgebra homomorphism of degree 0,
$f^*:\tilde T(M)\rightarrow W,$ such that $f^*(m)=f(m)$,  for all
$m\in M$.    The fact that $f$ is an Akivis homomorphism implies
that $I\subseteq {\rm ker}\, f^*$.    Hence, there is an
homomorphism of superalgebras of degree 0, $\tilde f:\tilde
U(M)\rightarrow W$ such that $\tilde f(\iota(m))=f^*(m)=f(m)$,  all
$m\in M$.    As $\tilde U(M)$ is generated by K and $\iota(M)$,
the unicity of $\tilde f$ follows.\halmos

\vspace{5mm}

\section  {Enveloping superalgebras of Lie and Malcev superalgebras}

In [4],  P\'{e}rez-Izquierdo studied enveloping algebras of Sabinin
algebras, showing that these generalize the classical  notions of
enveloping algebras for the particular cases of Lie algebras and
Malcev algebras.    In this section,  we will study the connections
between the classical definitions and the  definition of enveloping
superalgebras of Lie  and Malcev superalgebras considered as Akivis
superalgebras.

  Given a Lie  superalgebra $L$,  denote by
  $(\tilde U(L),  \iota)$ its enveloping superalgebra as an
  Akivis superalgebra and by $(U(L),  \sigma)$ its classical universal enveloping superalgebra.
   Clearly $\sigma$ is an Akivis homomorphism from $L$ to $U(L)^A$.
   So there is a unique homomorphism  of superalgebras $\tilde \sigma$ such that $\tilde \sigma \iota=\sigma$.
   As $U(L)$ is generated by K and $\sigma(L)$,   $\tilde\sigma$ is surjective.
    So $U(L)$ is an epimorphic  image of $\tilde U(L)$ by $\tilde \sigma$,   i.
    e.,
     $$U(L)\simeq \tilde U(L)/ {\rm ker}\, \tilde \sigma \,.$$

Suppose now that $N$ is a Malcev superalgebra.    Superizing the
theory exposed in [5],   we can naturally define the classical
enveloping superalgebra of $N$ as the superalgebra $(\tilde T(N)/
\tilde I,  \tilde \iota)$,  where $\tilde I$ is the graded ideal of
$\tilde T(N)$ generated by the homogeneous elements $$a\otimes b
-(-1)^{\alpha \beta}b\otimes a-[a, b],\,\, \, (a,  x,
y)+(-1)^{\alpha \gamma} (x,  a,  y),\, \, \, (x,  a, y)+(-1)^{\alpha
\xi} (x,  y, a),$$ all $a\in N_{\alpha},\,  b\in N_{\beta}, \, x\in
\tilde T(N)_{\gamma}, \, y\in \tilde T(N)_{\xi}$, and where $(x,  y,
z)=(xy)z-x(yz)$ denotes the usual associator. The map $\tilde \iota$
is the composition of the canonical injection from $N$ to $\tilde
T(N)$ with the canonical epimorphism $\mu:\tilde T(N)\rightarrow
\tilde T(N)/ \tilde I.   $

We will prove that $\tilde T(N)/ \tilde I$ is an epimorphic image of
$\tilde U(N)$. For this,  we show that $I\subseteq \tilde I$ and so
$\mu$  gives rise to the epimorphism $\tilde \mu:\tilde T(N)/
I\rightarrow \tilde T(N)/ \tilde I$ defined by $\tilde\mu (n+I)=
n+\tilde I$, all $n \in \tilde T(N)$.    To see that $I\subseteq
\tilde I$ notice that,  for all homogeneous elements
 $x\in N_{\gamma},  y\in N_{\xi},  z\in N_{\alpha}$,  we have that
$$
SJ(x,  y,  z)-3(xy)z= (-1)^{\alpha(\gamma+\xi)}[(z,  x,
y)+(-1)^{\alpha \gamma} (x,  z,  y)]-(-1)^{\alpha\xi}[(x,  z,
y)+(-1)^{\alpha \xi} (x,  y,z)]
$$
 is an element of $\tilde I$ and so $SJ(x,  y,  z)-3(xy)z+ SJ(y,  z, x)-3(yz)x\in
\tilde I$.    Therefore,
 $SJ(x,  y,  z)\in \tilde I$ and $(x,  y,  z)\in \tilde I$.
  As the generators of $I$ are in $\tilde I$, we can conclude that $I\subseteq \tilde I$.

\vspace{5mm}

\section {Speciality of Akivis superalgebras}

The canonical filtration of $\tilde T(M)$,   $\tilde T_0(M)\subseteq
\tilde T_1(M)\subseteq .   .   .   $,   where $\tilde T_0(M)=K$, and
$\tilde T_n(M)=\oplus_{i=0}^n\tilde T^i(M)$,   $n>0$,  gives rise to
the canonical filtration of $\tilde U(M), $ $\tilde U_0(M)\subseteq
\tilde U_1(M)\subseteq ... $,   where $\tilde U_n(M)=\tilde
T_n(M)+I$.

Associated with this filtration there is the $Z$-graded
superalgebra,

$$
gr\tilde U(M)=\oplus_{n\in Z}(gr\tilde U(M))_n,$$ where $ (gr\tilde
U(M))_n=0, $ if $n<0$,   $(gr\tilde U(M))_0=K$, $(gr\tilde U(M))_n=
\tilde U_n (M)/ \tilde U_{n-1}(M)$,  for  $n\geq 1, $ with
multiplication  given by
$$
(a+\tilde U_{i-1}(M))(b+\tilde U_{j-1}(M))=ab+\tilde U_{i+j-1}(M)
$$
for $a\in \tilde U_i(M)$ and $b\in \tilde U_j(M)$.    For simplicity
we identify $(gr\tilde U(M))_1$ with $\iota(M)$.

Now consider the classical tensor algebra  $T(M)$ of $M$,   that is
naturally a $Z$-graded associative superalgebra,
$$
T(M)=\oplus_{n\in Z} T^n(M)
$$
where $T^n(M)=0$ if $n<0$,   $T^0(M)=K$,   $T^n(M)= M\otimes
M\otimes. . . \otimes M$ ($n$ times) if $n>0$ .    Let $J$ be  the
ideal of $T(M)$ generated by the homogeneous elements  $x\otimes y
-(-1)^{\alpha\beta} y\otimes x$,   all  $x\in M_{\alpha},   y\in
M_{\beta}, \,  {\alpha,  \beta}\in Z_2$. The associative Z-graded
quotient superalgebra $S(M)=T(M)/J$ is called the supersymmetric
superalgebra of $M$. Note that as associative Z-graded algebra the
homogeneous spaces of $S(M)$ are $S^n(M)=T^n(M)+J$ and as  a
superalgebra we have $S(M)_\alpha=\oplus_{n\in Z}
(\oplus_{\alpha_1+.   .   .   + \alpha_n=\alpha}
(M_{\alpha_1}\otimes.   .   .   \otimes M_{\alpha_n}+J))$.    Since
the generators of $J$ lie in $T^2(M), $ we identify $S^0(M)$ with
$K$ and $S^1(M)$ with $M$.

We will now construct from $S(M)$ a nonassociative superalgebra
$V(M)$ which will play in this work  the role that the symmetric
algebra plays in the classical case of the PBW Theorem.

We define the $Z$-graded supervector space $V(M)=\oplus_{n\in Z}V^n
(M), $ where the subspaces $V^n (M)$  are defined by
 $$
V^n(M)= S^n(M), \, \, {\rm if}   \, \, n\leq 3,  \, \, {\rm and} \,
\, V^n(M)= \oplus_{i=1}^{n-1} (V^i(M)\otimes V^{n-i}(M)), \, \, {\rm
if}  \, \, n>3.
$$

We turn $V(M)$ into a superalgebra by defining the multiplication
for homogeneous elements $v_i\in V^i(M),  v_j\in V^j(M)$  by $ v_i.
v_j=v_iv_j$ if $i+j\le 3$ and $v_i.   v_j=v_i\otimes v_j$ if
${i+j}>3$ (where juxtaposition of elements  means the product of
these elements in  $S(M)$).

\begin{lemma}
 The superalgebra $V(M)$ is the enveloping Akivis
superalgebra of the trivial Akivis superalgebra with underlying
vector space $M$,   i.e.,  $V(M)\cong \tilde T(M)/I^*$ where  $I^*$
is the ideal of $\tilde T(M)$ generated by the homogeneous elements
$x\otimes y -(-1)^{\alpha\beta} y\otimes x,  \,  (x\otimes y)\otimes
z -x\otimes (y\otimes z)$,   all $x\in M_\alpha,   y\in M_\beta,
z\in M_\gamma$. \end{lemma}
 \noindent{\bf Proof:}  The inclusion map
$\pi:M\rightarrow V(M)^A$ is an Akivis homomorphism. Therefore, as
$\tilde T(M)/I^*$ is the enveloping superalgebra of the trivial
Akivis superalgebra obtained from $M$,   there is a superalgebra map
of degree 0, $\tilde \pi: \tilde T(M)/I^*\rightarrow V(M)$ such that
$\tilde \pi(m+I^*)= m$,   all $m\in M$.

Notice that in $\tilde T(M)/ I^*$ we have $x\otimes y
-(-1)^{\alpha\beta} y\otimes x+I^*=0$ and $(x\otimes y)\otimes z
-x\otimes (y\otimes z)+I^*=0$, for homogeneous elements $x, y, z$.
Hence,  using the universal property of the tensor products,   we
may define linear maps $\rho_1:M\rightarrow \tilde T(M)/ I^*$,
$\rho_2:M\otimes M\rightarrow \tilde T(M)/ I^*$,  $\rho_3:M\otimes
M\otimes M\rightarrow \tilde T(M)/ I^*$ by $\rho_i(a)=a+I^*,$  $i=1,
2,  3$.    Clearly,  $J\cap M\otimes M \subseteq {\rm ker}\, \rho_2$
and $J\cap M\otimes M\otimes M \subseteq {\rm ker}\, \rho_3$.
Therefore, there are linear maps $\tilde \rho_i:V^i(M)\rightarrow
\tilde T(M)/ I^*$ defined by $\tilde \rho_i(a+J)=a+I^*$,  i=1,  2,
3.    Using once more the universal property of tensor products and
induction, one can extend these maps to an homomorphism of
superalgebras of degree 0,  $\tilde \rho:V(M)\rightarrow \tilde
T(M)/ I^*$ satisfying $\tilde \rho(m)=m+I^*$ for all $m\in M$. The
maps $\tilde\rho$ and $\tilde \pi$ are inverse of each other. So
$\tilde \pi$ is an isomorphism.\halmos

\begin{lemma} There exists an epimorphism of $Z$-graded superalgebras
 $\,\tilde \tau: \tilde T(M)/I^*\rightarrow gr\tilde U(M)\,,$ of degree $0$,
such that $\tilde \tau(m+I^*)=\iota(m)$,  all $m\in M$.
\end{lemma}
\noindent{\bf Proof:}  Consider the natural epimorphism of
$Z$-graded algebras $\tau:\tilde T(M)\rightarrow gr\tilde U(M)$
given by $\tau(a)=(a+I)+\tilde U_{n-1}(M)$,   for each $a\in \tilde
T^n(M)$. Notice that,  since we identify $(gr\tilde U(M))_1$ with
$\iota(M)=M+I$,   then $\tau(m)=\iota(m)$, for all $m\in M$. To see
that $\tau$ preserves the gradation, recall that for $\alpha\in
Z_2$,  $\,\tilde T(M)_\alpha=\oplus_{n\geq 0} \tilde T^n(M)_\alpha$
and $$(gr\tilde U(M))_\alpha=\oplus_{n\geq 0} (\tilde U_n(M)/\tilde
U_{n-1}(M))_\alpha=\oplus _{n\geq 0}(\oplus_{i=0}^n(\tilde
T^i(M)_\alpha+I)+\tilde U_{n-1}(M)).$$ So if $a\in \tilde
T^n(M)_\alpha$,   then $\tau(a)=(a+I)+\tilde U_{n-1}(M)\in (\tilde
T^n(M)_\alpha+I)+\tilde U_{n-1}(M)$.    Hence $$\tau(\tilde
T(M)_\alpha)=\sum_{n\geq 0}\tau (\tilde T^n(M)_\alpha)\subseteq
\oplus_{n\geq 0}(\tilde T^n(M)_\alpha+I)+\tilde U_{n-1}(M))\subseteq
(gr\tilde U(M))_\alpha, $$   as desired. Now for any $x\in M_\alpha,
y\in M_\beta,  z\in M_\gamma$,   consider $\bar x=\iota(x),$ $ \bar
y=\iota(y),$  $ \bar z=\iota(z)$ $\in (gr\tilde U(M))_1$.    Then,
in $gr\tilde U(M)$ there holds
$$
\bar x \bar y-(-1)^{\alpha\beta} \bar y\bar x=
(\iota(x)\iota(y)-(-1)^{\alpha\beta} \iota(y)\iota(x))+ \tilde
U_1(M)=\iota([x,  y])+\tilde U_1(M)=0;
$$
$$
(\bar x\bar y)\bar z-\bar x(\bar y\bar z)= (\iota(x)\iota(y))\iota
(z)-\iota(x)(\iota(y)\iota(z))+\tilde U_2(M)=\iota(A(x,  y,
z))+\tilde U_2(M)=0.
$$

This implies that $\,I^*\subseteq {\rm Ker}\, \tau\,$.    Therefore,
there is an epimorphism of $Z$-graded superalgebras $\tilde
\tau:\tilde T(M)/I^*\rightarrow gr \tilde U(M)$ such that $\tilde
\tau (m+I^*)=\tau(m)=\iota(m)$,   for all $m\in M$.\halmos

\vskip5mm
 From the two previous lemmas,  we know that the composite map
 $\tilde \tau {\tilde\pi}^{-1}:V(M)\rightarrow gr\tilde U(M)$ is an
 epimorphism of superalgebras satisfying $\tilde \tau {\tilde\pi}^{-1}(m)=\iota(m)$ for all $m\in M$.
  It is our aim,  to prove that this epimorphism is in fact an isomorphism.    For this we need to
  endow $V(M)$ with a convenient superalgebra structure.
\vskip5mm

Let $\{e_r,  r\in \Delta\}$ be a basis of $M$ indexed by the totally
ordered set $\Delta=\Delta_0\cup \Delta_1$ satisfying the following:
$\{e_r:r\in \Delta_\alpha\}$ is a basis of $M_\alpha, \,  \alpha=0,
1$,  and $r<s$ if $r\in \Delta_0, \,  s\in \Delta_1$.    It is well
known that,  in these conditions,  $\{e_{r_1}e_{r_2}:r_1\le r_2$,
and $r_1<r_2$ if $r_1,  r_2\in \Delta_1\}$ is a basis of $V^2(M)$
and  $\{e_{r_1}e_{r_2}e_{r_3}:r_1\le r_2\le r_3$,  and $r_p<r_{p+1}$
if $r_p,  r_{p+1}\in \Delta_1\}$ is a basis of $V^3(M)$.
 In the supervector space $V(M)$ we define a new multiplication denoted by $*$ in the following way:
if $a\in V^i(M), \,  b\in V^j(M)$,  then $a*b=a\otimes b$ if
$i+j>3$; if $i+j\le 3$ the multiplication is defined on the basis
elements by the following identities (for simplicity we use $\bar r,
\bar s,  \bar k$ to denote the degrees of $e_r,  e_s,  e_k$,
respectively):

$$e_r*e_s=\cases {e_re_s,&if  $r\le s \, \,  {\rm and}\, \,  r\not= s \, \,  {\rm if}\, \,  r\in \Delta_1$;\cr
                  1/2[e_r,  e_r],& if $r=s \in \Delta_1$;\cr
                  (-1)^{\bar r\bar s} e_se_r +[e_r,  e_s],& if $r > s.$\cr}$$

 $$(e_re_s)*e_k= \cases { e_r e_s e_k,& if $r\le s\le k\, \,  {\rm and} \, \, k\not= s \,  \,{\rm if}\, \,  s\in \Delta_1 ;$\cr
                          A(e_r,  e_s,  e_s)+1/2 e_r*[e_s,  e_s],&  if  $r<s=k \, \,{\rm and} \, \,s\in \Delta_1;$\cr
                         (-1)^{\bar k \bar s} e_r e_k e_s +e_r*[e_s,  e_k]+\cr
+A(e_r,  e_s,  e_k)-(-1)^ {\bar s\bar k} A(e_r,  e_k,  e_s),&  if
$r\le k<s \, \,{\rm and}\, \, r\not=k\, \, {\rm if} \, \,r\in
\Delta_1;$\cr
                        -1/2[e_r,  e_r]*e_s+e_r*[e_s,  e_r]+\cr
+A(e_r,  e_s,  e_r)+A(e_r,  e_r,  e_s),& if $r=k<s \, \,{\rm and} \,
\,r\in \Delta_1;$\cr
                         (-1)^{\bar k(\bar r+\bar s}) e_ke_r e_s +(-1)^{\bar k \bar s}[e_r,  e_k]*e_s +\cr
+e_r*[e_s,  e_k]-(-1)^ {\bar s \bar k} A(e_r,  e_k,  e_s)+\cr +
A(e_r,  e_s,  e_k),& if  $k<r\le s.$   \cr}$$

$$e_r*(e_s e_k)=\cases {e_r e_s e_k-A(e_r,  e_s,  e_k),& if $r\le s\le k \, \,{\rm and} \, \,r\not=s \, \,{\rm if}\, \, s\in \Delta_1;$ \cr
                         1/2[e_r,  e_r]*e_k-A(e_r,  e_r,  e_k),&   if $r=s<k \, \,{\rm and}\, \, r\in \Delta_1;$\cr
                        (-1)^{\bar r \bar s} e_se_r e_k -A(e_r,  e_s,  e_k)+\cr +[e_r,  e_s]*e_k,& if $s<r\le k\, \, {\rm and} \, \,r\not= k \, \,{\rm if}\, \, r\in \Delta_1;$\cr
                         1/2 (-1)^{\bar r\bar s} e_s*[e_r,  e_r]+[e_r,  e_s]*e_r -\cr -A(e_r,  e_s,  e_r)+(-1)^{\bar r \bar s}A(e_s,  e_r,  e_r),& if $s<r=k \, \,{\rm and}\, \, r\in \Delta_1;$\cr
                        (-1)^{\bar r(\bar k+\bar s)} e_s e_k e_r +(-1)^{\bar r\bar s} e_s*[e_r,  e_k]+\cr
 +[e_r,  e_s]*e_k+(-1)^{\bar r\bar s}A(e_s,  e_r,  e_k)-\cr
-A(e_r,  e_s,  e_k) - (-1)^ {(\bar s+\bar k) \bar r}A(e_s,  e_k,
e_r),& if  $s\le k<r.$ \cr}$$
 Note that if we consider a basis element $e_pe_q$ of $V^2(M)$ we always assume $p\leq q$ and $p\not= q$ if $p, q\in \Delta_1$.

With this multiplication $V(M)$ becomes a superalgebra that  will be
denoted by $\tilde V(M)$.

\begin{lemma} There is an homomorphism  of superalgebras $\hat
\epsilon: \tilde U(M)\rightarrow \tilde V(M)$, of degree 0,
satisfying $\hat \epsilon(\iota (m))=m$,  for all $m\in M$.
\end{lemma}
\noindent{\bf Proof:} Denote the operations in the Akivis
superalgebra $\tilde V(M)^A$ by
$$
<x,  y>= x*y-(-1)^{\alpha\beta}y*x\,,\,\, \,\,\,\, <x,  y,  z>=
(x*y)*z-x*(y*z),
$$
for $x\in M_\alpha, \,  y\in M_\beta,  \,  z\in M_\gamma$.

We  start by proving that the inclusion map $\epsilon:
M\rightarrow\tilde V(M)^A$ is an homomorphism of Akivis
superalgebras.    For this,   it is enough to show that  $$\,[e_r,
e_s]=<e_r,  e_s>\,\,\,\, \mbox{and} \,\,\,\,A(e_r,  e_s,  e_k)=<e_r,
e_s, e_k>\,,\,$$ for the basis elements $e_r,  e_s,  e_k$ considered
above. It is quite simple to see that the first of these two
inequalities holds. For the second one,  we have to consider several
cases. Here we present only four of them,   being the other cases
similar.

\begin{enumerate}

\item  $r=s=k\in \Delta_1$:

$<e_r,  e_r,  e_r>=(e_r*e_r)*e_r-e_r*(e_r*e_r)=1/2([e_r,
e_r]*e_r-e_r*[e_r,  e_r]).$

As $[e_r,  e_r]\in M_0$,   we have $[e_r,  e_r]=\sum_{t\in
\Delta_0}\alpha_t e_t$ (sum with finite support),  for scalars
$\alpha_t \in K$.    Therefore,   as $t<r$ for any $t\in \Delta_0$,
there holds   $$[e_r,  e_r]*e_r=\sum_{t\in \Delta_0}\alpha_t e_t*e_r
= \sum_{t\in \Delta_0}\alpha_t e_t e_r.  $$ In a similar way,   we
see that
$$
e_r*[e_r,  e_r]= \sum_{t\in \Delta_0}\alpha_t(e_t e_r+[e_r,  e_t])=
\sum_{t\in \Delta_0}\alpha_t e_t e_r+ [e_r,  [e_r,  e_r]].
$$
Hence $<e_r,  e_r,  e_r>=-1/2[e_r,  [e_r,  e_r]]=1/2[[e_r,  e_r],
e_r].   $ On the other hand,
 from the definition of Akivis superalgebra,  we have that $SJ(e_r,  e_r,  e_r)=3[[e_r,  e_r],  e_r]=
 6 A(e_r,  e_r,  e_r)$.    Thus $<e_r,  e_r,  e_r>=A(e_r,  e_r,  e_r).   $

\item  $r\leq k<s$ and $r\not=k$ if $r\in \Delta_1$:
$$\begin{array}{l}
<e_r,  e_s,  e_k> = (e_r e_s)*e_k-e_r*((-1)^{\bar s \bar k}
e_ke_s+[e_s,  e_k])= (-1)^{\bar s \bar k}e_r e_k e_s+A(e_r,  e_s,
e_k)- \\-(-1)^{\bar s \bar k}A(e_r,   e_k,  e_s)+ e_r*[e_s,  e_k]
-e_r*[e_s, e_k]- (-1)^{\bar s \bar k}e_re_ke_s+ (-1)^{\bar s \bar k}
A(e_r,e_k,e_s)=\\
=A(e_r,  e_s,  e_k). \end{array}$$

\item   $r=k<s$ and $r\in \Delta_1$:
$$\begin{array}{l}
<e_r,  e_s,  e_r>=-1/2[e_r, e_r]*e_s+e_r*[e_s, e_r]+\\+A(e_r, e_s,
e_r)+A(e_r, e_r, e_s)+ e_r*(e_re_s)-e_r*[e_s, e_r]=\\= A(e_r, e_s,
e_r)+A(e_r, e_r, e_s)-1/2[e_r, e_r]*e_s+ 1/2[e_r, e_r]*e_s- A(e_r,
e_r, e_s)=\\=A(e_r, e_s, e_r). \end{array}$$

\item   $k<s<r$:
$$\begin{array}{l} <e_r,  e_s,  e_k> = (-1)^{\bar s \bar r}(e_s
e_r)*e_k+[e_r, e_s]*e_k- (-1)^{\bar s \bar k}e_r* (e_k
e_s)-e_r*[e_s,  e_k]=\\
=(-1)^{(\bar s+\bar k) \bar r}[e_s,  e_k]*e_r+(-1)^{\bar s \bar
r}e_s* [e_r,   e_k]+[e_r,  e_s]*e_k - e_r*[e_s,  e_k]-
\\-(-1)^{(\bar s+\bar r) \bar k} e_k*[e_r,  e_s]-
(-1)^{\bar s \bar k}[e_r,  e_k]*e_s+ (-1)^{\bar s \bar r}A(e_s, e_r,
e_k)- \\-(-1)^{(\bar s+\bar k) \bar r} A(e_s,  e_k,  e_r)+
(-1)^{{\bar s \bar k}+(\bar s+\bar k) \bar r} A(e_k,  e_s,
e_r)-(-1)^{(\bar s+\bar r) \bar k}A(e_k,  e_r,  e_s)+ \\+(-1)^{\bar
s \bar k} A(e_r, e_k, e_s)=\,\,\,\,(\, \, \, {\rm since} \, \, \,
[x, y]=<x, y>, \, \, \, {\rm all}\,
 \, \,  x, y\in M)
\\
=SJ(e_r,  e_s,  e_k)-((-1)^{\bar s \bar r} A(e_s,  e_r,
e_k)-(-1)^{(\bar s+\bar k) \bar r} A(e_s,  e_k,  e_r)+\\+(-1)^{\bar
s \bar k+(\bar s+ \bar k) \bar r}A(e_k,  e_s,  e_r) -(-1)^{(\bar
s+\bar r) \bar k}A(e_k,  e_r,  e_s) +(-1)^{\bar s \bar k} A(e_r,
e_k, e_s)) = \\=A(e_r,  e_s,  e_k)\,\,\,\, (\, {\rm by\, \, the\, \,
definition \, \,of\, \, Akivis \, \,superalgebra}).
\end{array}$$
 \end{enumerate}

As $\epsilon$ is an  Akivis homomorphism,  from the definition of
enveloping superalgebra,  there is an  homomorphism of superalgebras
of degree 0, $\hat \epsilon:\tilde U(M)\rightarrow \tilde V(M)$
satisfying $\hat \epsilon(\iota (m))=m$,   all $m\in M.$ \halmos

\begin{theorem} The $Z$-graded superalgebras $V(M)$ and $gr \tilde
U(M)$ are isomorphic. \end{theorem}

\noindent {\bf Proof:}  The superalgebra $\tilde V(M)$ has a natural
filtration defined by the sequence of subspaces $\tilde
V_n(M)=\oplus_{i=0}^n V^i(M)$. So we may consider the associated
$Z$-graded superalgebra $gr\tilde V(M)$. ( As usual we identify
$(gr\tilde V(M))_1$ with $M$). Since the map $\hat\epsilon$,
considered in the previous lemma, is an homomorphism of
superalgebras, we have $\hat\epsilon(\tilde U_n(M))\subseteq \tilde
V_n(M)$ . So, we may define $\tilde \epsilon:gr \tilde
U(M)\rightarrow gr\tilde V(M)$ by $\tilde\epsilon(a_i+\tilde
U_{i-1}(M))= \hat \epsilon (a_i)+\tilde V_{i-1}(M)$. This map is an
homomorphism of  $Z$-graded superalgebras of degree 0 and satisfies
$\tilde\epsilon (\iota(m))=m,$ all $m\in M$.

We now return to the algebra $V(M)$.   For $n\geq 1$, the map
$\mu_n:V^n(M)\rightarrow (gr\tilde V(M))_n$ defined by
$\mu_n(v)=v+\tilde V_{n-1}(M)$ is an isomorphism of vector spaces.
So taking $\mu_0=Id_K$,   $\mu=\oplus_{n\geq 0} \mu_n
:V(M)\rightarrow gr\tilde V(M)$ is an isomorphism of $Z$- graded
vector spaces.    Looking at the formulas which define the
multiplication in $\tilde V(M)$,   it is easy to see that this is in
fact an isomorphism of $Z$-graded superalgebras.    The composite
homomorphisms $\mu^{-1}\tilde \epsilon: gr \tilde U(M)\rightarrow gr
\tilde V(M)\rightarrow V(M)$ and  $\tilde\tau{\tilde \pi}^{-1}: V(M)
\rightarrow \tilde T(M)/I^*\rightarrow gr \tilde U(M)$ (recall the
preceeding lemmas) are inverse of each other. So the  result
follows.\halmos \vspace{5mm}

The following results are immediate consequence of this theorem.

\begin{corollary} The canonical map $\iota:M\rightarrow \tilde U(M)$ is
injective. \end{corollary}

\begin{corollary} Any Akivis superalgebra, defined over a field of
characteristic different from 2 and 3, is special. \end{corollary}

\vskip5mm

\section  {BIBLIOGRAPHY}

[1] Akivis M.   A.,    "Local algebras of a multidimensional three
web" (Russian), {\it Sibirsk.    Mat.   Zh.   } {\bf 17 }(1),
(1976) 5\--11;

[2]  Albuquerque H.    and Shahn Majid,  "Quasialgebra structure of
the octonions",   {\it Journal of Algebra} {\bf 220},   (1999)
188\--224;

[3] Albuquerque H.,   A.    Elduque and Jos\'{e} P\'{e}rez-Izquierdo,
"$Z_2$- quasialgebras", {\it Com.    in Algebra} {\bf 30 } (5),
(2002)  2161\--2174;

[4] P\'{e}rez-Izquierdo J.,   "Algebras,   hyperalgebras, nonassociative
bialgebras and loops", {\it Advances in Mathematics} {\bf 208}
(2007) 834\--876;

[5] P\'{e}rez-Izquierdo J.    and Ivan Shestakov,   "An envelope for
Malcev algebras", {\it J.    of Algebra} {\bf 272} (1) (2005)
379-393;

[6] Shestakov I.,    "Every Akivis Algebra is Linear", {\it
Geometriae Dedicata} {\bf 77} (2) (1999) 215\--223;

[7] Shestakov I.   and U.   U.   Umirbaev,  "Free Akivis Algebras,
primitive elements and hyperalgebras", {\it J.   Algebra} {\bf
250}(2) (2002) 533-548.

[8] Sheunert M.,    " The theory of Lie superalgebras", {\it LNM}
{\bf 716} Springer Verlag,   Berlin 1979.

\end{document}